\documentclass[a4paper,10pt]{article}
\usepackage[T1]{fontenc}
\usepackage[french]{babel}
\usepackage{amssymb}
\usepackage{amsmath}
\usepackage{theorem}

\newtheorem{atheorem}{\hspace{\parindent}Th\'{e}or\`{e}me}

\newcommand{\sophistnorm}[2]{\left\Vert\,#1\,\right\Vert_{#2}}
\newcommand{\maxnorm}[1]{\left\Vert\,#1\,\right\Vert_{\infty}}

\newcommand{\scal}[1]{\langle\,#1\,\rangle}

\numberwithin{equation}{section} \theoremstyle{plain} \theoremheaderfont{\scshape}
\theorembodyfont{\itshape}

\newcommand{\eps}{\varepsilon}
\newcommand{\la}{\lambda}

\newcommand{\de}{\Delta}

\newcommand{\spl}{\mathbb{S}\,(L,\de)}
\newcommand{\pr}{P_\mathbb{S}}
\newcommand{\prj}{P_{\mathbb{S}\,(L,\de)}}
\DeclareMathOperator{\prs}{P_\mathbb{S}}

\newcommand{\lf}{\lambda_1}
\newcommand{\ls}{\lambda_2}

\def\define{\mathrel{\stackrel{\rm def}=}}

\renewcommand{\phi}{\varphi}

\begin{document}

\renewcommand{\labelitemi}{\textbullet}

\thispagestyle{empty}
\begin{large}
\begin{center}

{Minist\`{e}re de l'\'{e}ducation de la F\'{e}d\'{e}ration de Russie}

{Universit\'{e} d'\'{e}tat de l'Oural}

\bigskip
{Facult\'{e} de math\'{e}matiques et m\'{e}canique}

{Chair d'analyse math\'{e}matique et th\'{e}orie des fonctions}

\vspace*{\stretch{1}}

{\Large\uppercase{La meilleure approximation moindres carr\'{e}s\\par les $L$-splines de
deuxi\`{e}me ordre}}

\vspace{4cm}

\begin{tabular}{l r}
\hspace*{-0.5cm}\parbox{0.3\linewidth}
{\vspace*{-58pt}Chef de la chair\\ professeur Ar\'{e}stov V.V.\\
\vspace{1cm}\\ \rule{5cm}{0.5pt}}

&

\hspace*{1cm}\begin{minipage}{0.7\linewidth}

Th\`{e}se pour l'obtention du magist\`{e}re\\
sp\'{e}cialit\'{e} \guillemotleft\thinspace Math\'{e}matiques, math\'{e}matiques
appliqu\'{e}es\thinspace \guillemotright \\
pr\'{e}sent\'{e}e et soutenue par\\
Kaioumov Alexandre Rachidovitch

\bigskip
\bigskip
\bigskip
Directeur de th\`{e}se ---\\
membre correspondant de l'Acad\'{e}mie des Science de la Russie,\\
professeur, docteur \`{e}s sciences physico-math\'{e}matiques,\\
chef du Departement d'approximation\\
de l'Institut de mathematiques et mecanique\\
de l'Academie des science de la Russie\\
Soubbotine Youri Nikolai\'{e}vitch
\end{minipage}
\end{tabular}

\end{center}

\vspace*{\stretch{2}}

\begin{center}
\'{E}kat\'{e}rinebourg, 2003
\end{center}
\end{large}

\newpage
\thispagestyle{empty}

\begin{center}
{\Large\bf{R\'{e}sum\'{e}}}
\end{center}
\bigskip

{\bf Mots-clefs:} {\uppercase{$L$-spline, spline exponentielle, spline
trigonom\'{e}trique, meilleure approximation moindres carr\'{e}s, projecteur orthogonal,
constante de Lebesgue.}}

\bigskip
{\bf Sujet de la th\`{e}se.} Dans cette th\`{e}se nous \'{e}tudions l'op\'{e}rateur de
l'approximation moindres carr\'{e}s par les $L$-splines g\'{e}n\'{e}r\'{e}es par un
op\'{e}rateur diff\'{e}rentiel lin\'{e}aire de deuxi\`{e}me ordre --- autrement dit,
l'op\'{e}rateur de projection orthogonale sur l'espace de ces $L$-splines. Nous
consid\'{e}rons comme l'op\'{e}rateur diff\'{e}rentiel les op\'{e}rateurs
g\'{e}n\'{e}rant les $L$-splines exponentielles et trigonom\'{e}triques. Nous traitons la
question suivante: la norme uniforme d'un tel op\'{e}rateur de projection, est-elle
born\'{e}e uniform\'{e}ment en tous les partitionnements du segment?

\vspace{4pt} {\bf R\'{e}sultats principaux.} Nous d\'{e}montrons que les op\'{e}rateurs
de projection orthogonale sur de diff\'{e}rents espaces de $L$-splines sont born\'{e}s
uniform\'{e}ment en tous les partitionnements du segment. Pour les splines exponentielles
avec les exposants d'une m\^{e}me valeur absolue, une estimation exacte est obtenue pour
la norme du projecteur. Pour les splines exponentielles avec les exposants arbitraires de
signes diff\'{e}rents, nous d\'{e}montrons que le projecteur est uniform\'{e}ment
born\'{e}. Pour les splines trigonom\'{e}triques avec des contraintes habituelles sur le
segment, une estimation est obtenue pour la norme du projecteur. Nous traitons aussi un
exemple mod\`{e}le dans le but de montrer que ces contraintes sur le segment ne sont pas
du tout n\'{e}cessaires pour que le projecteur orthogonal soit born\'{e}.

\vspace{4pt} {\bf M\'{e}thodes employ\'{e}es.} Nous examinons les matrices de Gram des
bases de $LB$-splines, aussi que leurs inverses. Dans le cas des splines
trigonom\'{e}triques, au cours de l'\'{e}tude de l'exemple mod\`{e}le, nous
\'{e}tablissons un ph\'{e}nom\`{e}ne d'interf\'{e}rence, qui permet \`{a} l'approximant
trigonom\'{e}trique de rester born\'{e}, tandis que ses coefficients et les splines de
base tendent vers l'infini.

\vspace{4pt} {\bf Histoire du probl\`{e}me et port\'{e}e des r\'{e}sultats.} L'\'{e}tude
des normes des projecteurs orthogonaux sur de diff\'{e}rents espaces fonctionnels est un
probl\`{e}me classique de la th\'{e}orie d'approximation bien important en lui-m\^{e}me.
Appliqu\'{e} aux splines polynomiales, ce probl\`{e}me est connu sous le nom de {\it
probl\`{e}me de Carl de Boor}.

Dans le cas des splines, la norme du projecteur orthogonal sur l'espace des splines de
degr\'{e} $k-1$ est li\'{e}e \`{a} la norme de l'op\'{e}rateur d'interpolation par les
splines de degr\'{e} $2k-1$, ce qui accentue l'importance de son \'{e}tude. Professeur
C.~de~Boor a lanc\'{e} une \'{e}tude syst\'{e}matique de ce probl\`{e}me, qui l'a
amen\'{e} \`{a} la conclusion que la norme du projecteur devait \^{e}tre born\'{e}e
ind\'{e}pendamment du partitionnement du segment, une conjecture qu'il a formul\'{e}e en
1972.

Malgr\'{e} de nombreux r\'{e}sultats partiels, le probl\`{e}me de C.~de~Boor restait
ouvert pendant plus de 25 ans, jusqu'\`{a} ce qu'il n'ait \'{e}t\'{e} r\'{e}solu
positivement par Dr. A.Yu. Chadrine en 1999.

Le probl\`{e}me consider\'{e} dans cette th\`{e}se constitut donc la premi\`{e}re
\'{e}tape dans la g\'{e}n\'{e}ralisation du probl\`{e}me de C.~de~Boor aux $L$-splines.

\newpage

\section{Introduction}
\subsection{Position du probl\`{e}me}

Nous introduisons les notations suivantes. Fixons un segment $[a,b]$ et d\'{e}signons par

\begin{itemize}
\item $\de=\{a=t_1<t_2<\dots<t_{n-1}<t_n=b\}$ --- un partitionnement du segment $[a,b]$,

\item $\de_i=[t_i,t_{i+1}]\ (i=1,\dots,n-1)$ --- les sous-segments de ce partitionnement,

\item $L=D^k+a_{k-1}D^{k-1}+\dots+a_1D+a_0I,\ a_i\in\mathbb{R}$ --- un op\'{e}rateur
diff\'{e}rentiel lin\'{e}aire d'ordre $k$ \`{a} coefficients constants,

\item $N_L= \left\{ f\in C^k[a,b] : Lf(x)=0,\ x\in[a,b] \,\right\}$ --- le noyau de
l'op\'{e}rateur $L$.
\end{itemize}

Alors {\it l'espace des $L$-splines g\'{e}n\'{e}r\'{e} par le partitionnement $\de$ et
l'op\'{e}rateur diff\'{e}rentiel $L$} est d\'{e}fini par
\begin{equation*}
\spl = \left\{ s\in C^{k-2}[a,b] : s|_{\de_i}\in N_L,\ i=1,\dots,n-1 \,\right\}.
\end{equation*}

Consid\'{e}rons l'op\'{e}rateur $\pr=\prj$ de projection orthogonale (relative au produit
int\'{e}rieur dans l'espace $L_2[a,b]$) de l'espace $C[a,b]$ sur son sous-espace $\spl$:
\begin{equation*}
\pr : C[a,b]\ni f \longmapsto \prs{f}\in\spl,
\end{equation*}
o\`{u} $\prs{f}$ est uniquement d\'{e}fini par l'une quelconque des deux conditions
suivantes:
\begin{enumerate}
\item $\scal{f,s}=\scal{\prs{f},s},\quad\forall s\in\spl$,

\item $\sophistnorm{f-\prs{f}}{2} = \min\limits_{s\in\spl}\sophistnorm{f-s}{2}$.
\end{enumerate}
Nous nous occupons de la norme de l'op\'{e}rateur $\pr$ \`{a} titre d'un op\'{e}rateur de
$C[a,b]$ \`{a} valeurs dans $C[a,b]$:
\begin{equation*}
\maxnorm{\pr} = \sup\big\{\, \maxnorm{\prs{f}} : f\in C[a,b],\ \maxnorm{f}\leqslant 1
\,\big\}.
\end{equation*}
Il est ais\'{e} de voir que lorsque le partitionnement $\de$ est fix\'{e},
l'op\'{e}rateur $\prj$ est born\'{e}, c'est-\`{a}-dire a une norme finie. Nous sommes
donc naturellement amen\'{e} \`{a} la question suivante: la famille d'op\'{e}rateurs de
projection $\prj$, est-elle born\'{e}e {\it uniformement en tous les partitionnements
$\de$ du segment $[a,b]$}\,?

\subsection{Histoire du probl\`{e}me}

Le probl\`{e}me que nous venons de formuler est la g\'{e}n\'{e}ralisation au cas des
$L$-splines du probl\`{e}me correspondant pour les splines polynomiales: on obtient le
dernier en posant $L\equiv D^k$ dans le premier (les $L$-splines ainsi obtenues
co\"{i}ncideraient avec les splines polynomiales d'ordre $k$, de degr\'{e} $k-1$).

Faisons un expos\'{e} de l'histoire du probl\`{e}me pour les splines polynomiales.
L'hypoth\`{e}se que les normes uniformes des projecteurs orthogonaux
$P_{\mathbb{S}_k(\de)}$ sur les espaces $\mathbb{S}_k(\de)$ des splines polynomiales sont
born\'{e}es ind\'{e}pendamment des partitionnements $\de$ du segment $[a,b]$ a
\'{e}t\'{e} \'{e}nonc\'{e}e par Carl de Boor dans \cite{deBoor-3} en 1972. C'est pourquoi
cette hypoth\`{e}se \'{e}tait connue sous le nom de {\it probl\`{e}me de C.~de~Boor}.
D\'{e}signons
\begin{equation*}
c_k = \sup_{\de}\maxnorm{P_{\mathbb{S}_k(\de)}}.
\end{equation*}

Dans le cas de $k=2$ (splines lin\'{e}aires ou lignes bris\'{e}es), Z.~Ciezielsky
\cite{Ciez} a montr\'{e} en 1963 l'estimation $c_2 \leqslant 3$. K.I.~Oskolkov \cite{Osk}
a obtenu en 1979 un r\'{e}sultat d'o\`{u} il d\'{e}coule que $c_2=3$.

Dans le cas de $k=3$ (splines paraboliques) C.~de~Boor \cite{deBoor-1} a montr\'{e} en
1968 l'estimation $c_3 \leqslant 30$.

Dans le cas de $k=4$ (splines cubiques) C.~de~Boor \cite{deBoor-2} a annonc\'{e} en 1979,
en partant de r\'{e}sultats d'exp\'{e}riences num\'{e}riques, l'assertion que $c_4
\leqslant \frac{245}{3}=81\frac{2}{3}$. Autant que je sache, une d\'{e}monstration exacte
de cette assertion n'a jamais \'{e}t\'{e} publi\'{e}e.

Apr\`{e}s cela, il n'y avait eu aucun progr\`{e}s pendant 20 ans --- jusqu'au 1999, quand
le probl\`{e}me de C.~de~Boor a \'{e}t\'{e} resolu par Dr. A.Yu. Chadrine dans
\cite{Shad}. Il a verifi\'{e} que $c_k\!<\!+\infty\ (k \geqslant 2)$, bien qu'il n'est
pas arriv\'{e} \`{a} obtenir une estimation sup\'{e}rieure constructive. En utilisant la
m\'{e}thode d'Oskolkov, il a n\'{e}anmoins obtenu une estimation inf\'{e}rieure:
$c_k\geqslant 2k-1$. En partant de r\'{e}sultats d'exp\'{e}riences num\'{e}riques et de
quelques consid\'{e}rations th\'{e}oriques, il a \'{e}nonc\'{e} la conjecture que $c_k =
2k-1\ (k \geqslant 2)$.

\medskip

Autant que je sache, le probl\`{e}me correspondant pour les $L$-splines n'a jamais
\'{e}t\'{e} examin\'{e}.

\subsection{Formulation des r\'{e}sultats}

Dans cette th\`{e}se nous nous limitons \`{a} l'analyse de $L$-splines de deuxi\`{e}me
ordre de type sp\'{e}cial, \`{a} savoir des splines {\it exponentielles} et {\it
trigonom\'{e}triques} de deuxi\`{e}me ordre. Quoique celles splines-ci soient analogiques
aux splines polynomiales lin\'{e}aires (lignes bris\'{e}es) qui ne pr\'{e}sentent aucune
difficult\'{e}, le passage des splines polynomiales aux $L$-splines rend la situation
assez complexe.

Les classes suivantes de $L$-splines sont trait\'{e}es et les th\'{e}or\`{e}mes suivants
sont d\'{e}montr\'{e}s.

\begin{enumerate}
\item {\it Splines exponentielles avec les exposants d'une m\^{e}me valeur absolue}.
Consid\'{e}rons les $L$-splines g\'{e}n\'{e}r\'{e}es par l'op\'{e}rateur diff\'{e}rentiel
$L=D^2-\la^2$. Une telle spline co\"{i}ncide sur chaque sous-segment du partitionnement
avec une combinaison lin\'{e}aire des fonctions $e^{+\la x}$ et $e^{-\la x}$.

\begin{atheorem} Pour les $L$-splines g\'{e}n\'{e}r\'{e}es par l'op\'{e}rateur diff\'{e}rentiel
$L=D^2-\la^2$ on a:
\begin{equation*}
c(+\la,-\la)\define \sup_{\de}\maxnorm{\prj}=3.
\end{equation*}
\end{atheorem}

Notons que la m\'{e}thode employ\'{e}e pour montrer l'exactitude de l'estimation
ci-dessus nous fournit une autre d\'{e}monstration du r\'{e}sultat d'Oskolkov (que
l'estimation par Ciezielsky pour les splines polynomiales lin\'{e}aires est exacte).

\medskip

\item {\it Splines exponentielles avec les exposants de signes diff\'{e}rents et de
diff\'{e}rentes valeurs absolues}. Consid\'{e}rons les $L$-splines g\'{e}n\'{e}r\'{e}es
par l'op\'{e}rateur diff\'{e}rentiel $L=(D-\lf)(D-\ls)$, o\`{u} $\ls<0,\ \lf>0$. Une
telle spline co\"{i}ncide sur chaque sous-segment du partitionnement avec une combinaison
lin\'{e}aire des fonctions $e^{\lf x}$ et $e^{\ls x}$.

\begin{atheorem} Pour les $L$-splines g\'{e}n\'{e}r\'{e}es par l'op\'{e}rateur diff\'{e}rentiel
$L=(D-\ls)(D-\lf)$, o\`{u} $\ls<0,\ \lf>0$, on a: il existe une constante $C(\lf,\ls)$
qui ne d\'{e}pend qu'au nombres $\lf,\ls$ et ne d\'{e}pend pas au segment $[a,b]$, telle
que la famille d'op\'{e}rateur de projection $\prj$ est born\'{e}e par la constante
$C(\lf,\ls)$ uniform\'{e}ment en tous les partitionnements $\de$ du segment $[a,b]$:
\begin{equation*}
c(\lf,\ls,a,b)\define \sup_{\de}\maxnorm{\prj} \leqslant C(\lf,\ls).
\end{equation*}
\end{atheorem}

\medskip

\item {\it Splines trigonom\'{e}triques}. Consid\'{e}rons les $L$-splines
g\'{e}n\'{e}r\'{e}es par l'op\'{e}rateur diff\'{e}rentiel $L=D^2+\la^2$. Une telle spline
co\"{i}ncide sur chaque sous-segment du partitionnement avec une combinaison lin\'{e}aire
des fonctions $\sin \la x$ et $\cos \la x$.

\medskip

Ce cas-ci diff\`{e}re des deux cas pr\'{e}c\'{e}dents et du cas des splines polynomiales:
comme le probl\`{e}me d'interpolation lagrangienne par les fonctions trigonom\'{e}triques
devient insoluble sur un segment d'une longueur de $\pi/\la$, les propri\'{e}t\'{e}s des
$L$-splines de base empirent lorsque le diam\`{e}tre du partitionnement tend \`{a} cette
valeur. C'est pourquoi en traitant les splines trigonom\'{e}triques on essaie
habituellement d'\'{e}viter la proximit\'{e} du diam\`{e}tre du partitionnement \`{a} la
valeur $\pi/\la$. De plus, on introduit souvent une contrainte encore plus forte:
$d\,(\de)\leqslant \frac{\pi}{2\la}$ (o\`{u} $d\,(\de)=\max\limits_{i=1,\dots,n-1}
(t_{i+1}-t_i)$
--- est le diam\`{e}tre du partitionnement $\de$). (Voir, par exemple, \cite{Schu}, page
455, o\`{u} \cite{Rozh}, page 135.)

\medskip

D\'{e}signons
\begin{equation*}
M\,(\la,\tau) \define \sup_{\la \cdot d(\de) < \tau} \maxnorm{\prj}.
\end{equation*}

\begin{atheorem} Pour les $L$-splines g\'{e}n\'{e}r\'{e}es par l'op\'{e}rateur diff\'{e}rentiel
$L=D^2+\la^2$ on a:
\begin{equation*}
M\,(\la,\tau) \leqslant \begin{cases} 2\,\dfrac{(1-\cos\tau)^2} {(\tau-\sin\tau)
\sin\tau}\,, & \text{{\it si} $0 < \tau \leqslant \dfrac{\pi}{2},$}
\\ 2\,\dfrac{(1-\cos\tau)^2} {(\tau-\sin\tau) \sin^2\tau}\,, &
\text{{\it si} $\dfrac{\pi}{2} < \tau < \pi$.} \label{est-first}
\end{cases}
\end{equation*}
\end{atheorem}

Mais quand $\tau \to \pi-0$ le deuxi\`{e}me membre de l'estimation ci-dessus est
\'{e}quivalent \`{a}
\begin{equation*}
\dfrac{8}{\pi}\dfrac{1}{(\pi-\tau)^2}.
\end{equation*}
Une question se pose d'une fa\c{c}on naturelle: est-ce que cette asymptotique correspond
\`{a} la conduite actuelle de la valeur $M\,(\la,\tau)$\,? En g\'{e}n\'{e}ral cette
question semble difficile, et nous consid\'{e}rons ici un exemple mod\`{e}le, quand parmi
tous les sous-segments il n'y a qu'un d'une longueur proche \`{a} $\frac{\pi}{\la}$.

\begin{atheorem}
Soit $[a,b]=[0,\frac{\pi}{\la}]$, un des sous-segments ait une longueur proche \`{a}
$\frac{\pi}{\la}$, la somme des longueurs de tous les autres sous-segments soit $\eps$.
Alors
\begin{equation*}
\max\limits_{i=1,\dots,n} \left|\, \pr{f}(t_i) \,\right| \leqslant \left(
\dfrac{38}{\pi}+O(\eps) \right)\maxnorm{f},
\end{equation*}
\end{atheorem}

\begin{atheorem}
Soit $[a,b]=[0,\frac{\pi}{\la}]$. Alors
\begin{equation*}
M\,(\la,\pi) < +\infty
\end{equation*}
\end{atheorem}

\end{enumerate}

Les r\'{e}sultats ci-dessus ont \'{e}t\'{e} publi\'{e}s dans \cite{Ka-1, Ka-2, Ka-3,
Ka-4}.


\begin{thebibliography}{00}

\bibitem{deBoor-1} {\sc C.~de Boor} On the convergence of odd-degree spline interpolation. --- Journal of Approximation
Theory, 1968, {\bf 1}, pp.~452--263.
\bibitem{deBoor-2} {\sc C.~de Boor} On a max-norm bound for the
least-square spline approximant. --- "Approximation and Function Spaces"\ (Z.~Ciezielsky,
Ed.), pp.~163--175, Proceedings of the International Conference (Gdansk, August 27--31,
1979), PWN, Warzawa, 1981.
\bibitem{deBoor-3} {\sc C.~de Boor} The quasi-interpolant as a tool
in elementary polynomial theory. --- "Approximatoin Theory"\ (G.G.~Lorentz, Ed.),
pp.~269--276, Academic Press, New York, 1973.
\bibitem{Ciez} {\sc Z.~Ciezielsky} Properties of the orthonormal
Franklin systems. --- Studia Mathematica, 1963, {\bf 23}, fasc.~2 pp.~141--157.
\bibitem{Osk} {\sc K.I.~Oskolkov} The upper bound of the norms of
orthogonal projections onto subspaces of polygonals. --- "Approximation Theory"\
(Z.~Ciezielsky, Ed.), Banach Center Publications, vol.~4, pp.~177--183, PWN -- Polish
Scientific Publishers, Warsaw, 1979.
\bibitem{Schu} {\sc L.L.~Schumaker} {\it Spline Functions: Basic
Theory} --- Wiley, New York, 1981.
\bibitem{Shad} {\sc A.Yu.~Shadrin} The $L_{\infty}$-norm of the
$L_2$-spline prjector is bounded independently of the knot sequence: A proof of de Boor's
conjecture. --- Acta Mathematica, 2001, {\bf 187}, pp.~59--137.
\bibitem{Rozh} {\sc A.I.~Rojenko} {\it Th\'{e}orie abstracte des splines. Manuel universitaire} ---
Novosibirsk, Russie, 1999 {\it (en russe)}.
\bibitem{Ka-1} {\sc A.R.~Kaioumov} La meilleure approximation moindres carr\'{e}s par
les L-splines de deuxi\`{e}me ordre. --- Travaux de l'\'{E}cole internationale
S.B.~Stetchkine sur la th\'{e}orie des fonctions (Russie, la r\'{e}gion de Tcheliabinsk,
Miass, 24 juillet -- 3 ao\^{u}t 1999), pp.~105--127, \'{E}kat\'{e}rinebourg, Russie, 1999
{\it (en russe)}.
\bibitem{Ka-2} {\sc A.R.~Kaioumov} La meilleure approximation moindres carr\'{e}s par
les L-splines de deuxi\`{e}me ordre. --- "Th\'{e}orie d'approximation de fonctions et
d'op\'{e}rateurs", pp.~82--83, Travaux de la Conf\'{e}rence internationale d\'{e}di\'{e}e
\`{a} la 80-e anniversaire de S.B.~Stetchkine (Russie, \'{E}kat\'{e}rinebourg, 28
f\'{e}vrier -- 3 mars 2000), \'{E}kat\'{e}rinebourg, Russie, 2000 {\it (en russe)}.
\bibitem{Ka-3} {\sc A.R.~Kaioumov} Une estimation d'ordre exact dans le probl\`{e}me de la
meilleure approximation moindres carr\'{e}s par les L-splines de deuxi\`{e}me ordre de
type sp\'{e}cial. --- Travaux du Centre math\'{e}matique N.I.~Lobatchevski, vol.~12,
pp.~93--94, "Lectures Lobatchevski -- 2001", Mat\'{e}riaux de l'\'{E}cole--conf\'{e}rence
internationale de jeunes chercheurs (Russie, Kazan, 28 novembre -- 1 d\'{e}cembre 2001),
Kazan, Russie, 2001 {\it (en russe)}.
\bibitem{Ka-4} {\sc A.R.~Kaioumov} Le probl\`{e}me de C.~de~Boor pour les L-splines de
deuxi\`{e}me ordre. --- "Probl\`{e}me des math\'{e}matiques pures et appliqu\'{e}es",
pp.~66--70, Travaux de la 33-e \'{E}cole--conf\'{e}rence r\'{e}gionale de jeunes
chercheurs (Russie, \'{E}kat\'{e}rinebourg, 28 janvier -- 1 f\'{e}vrier 2002),
\'{E}kat\'{e}rinebourg, Russie, 2002 {\it (en russe)}.

\end{thebibliography}
\end{document}